\documentclass[11pt,leqno]{amsart}
\usepackage{amsmath}
\usepackage{amsthm}
\usepackage{amssymb}
\theoremstyle{definition}

\theoremstyle{plain}

\begin{document}
\setcounter{page}{1}
\begin{flushleft}
\scriptsize{{Inequality Theory and Applications\/} {\bf 6}\,(2007), 
\\
Pages ???--???}
\end{flushleft}
\vspace{16mm}

\begin{center}
{\normalsize\bf On certain analytic functions of bounded boundary rotation} \\[12mm]
     {\sc K. O. BABALOLA$^{1,2}$} \\ [8mm]

\begin{minipage}{123mm}
{\small {\sc Abstract.}
   In this paper we introduce certain analytic functions of boundary rotation bounded by $k\pi$ which are of Caratheodory origin. With them we study two classes of analytic and univalent functions in  the unit disk $E=\{z\in \mathbb{C}\colon |z|<1\}$, which are also of bounded boundary rotation.}
\end{minipage}
\end{center}

 \renewcommand{\thefootnote}{}
 \footnotetext{$^1$\,Current Address: Centre for Advanced Studies in Mathematics,
Lahore University of Management Sciences, Lahore, Pakistan. E-mail:
kobabalola\symbol{64}lums.edu.pk}
 \footnotetext{$^2$\,Permanent Address: Department of Mathematics, University of Ilorin,
 Ilorin, Nigeria. E-mail:
kobabalola\symbol{64}gmail.com; babalola.ko\symbol{64}unilorin.edu.ng}

\def\iff{if and only if }
\def\S{Smarandache }
\newcommand{\norm}[1]{\left\Vert#1\right\Vert}

\vskip 12mm
{\bf 1. Introduction }
\medskip

Let $P(\beta)$ denote the Caratheodory family of functions:
$$p(z)=1+c_1z+c_2z^2+\cdots\eqno{(1.1)}$$
which are analytic in $E$ and satisfy Re $p(z)>\beta$, $0\leq\beta<1$, $z\in E$. For $\beta=0$, we write $P$ in place of $P(0)$. The function $L_{0,\beta}(z)=\beta+(1-\beta)L_0(z)$ plays the extremal role in $P(\beta)$ as does $L_0(z)=(1+z)/(1-z)$ in $P$. Functions in $P(\beta)$ have Herglotz representations
$$\int_0^{2\pi}\frac{1+(1-2\beta)ze^{-is}}{1-ze^{-is}}d\mu(s)$$
for $d\mu(s)\geq 0$ and $\int_0^{2\pi}d\mu(s)=1$. Denote by $M_k$, $k\geq 2$, the class of real-valued
functions $m(s)$ of bounded variation on $[0, 2\pi]$ which satisfy the conditions:
$$\int_0^{2\pi}dm(s)=2,\;\;\;\int_0^{2\pi}|dm(s)|\leq k.$$

Clearly $M_2$ is the class of nondecreasing functions on $[0,
2\pi]$ which satisfy $\int_0^{2\pi}dm(s)=2$.\vskip 2mm

If $m\in M_k$, with $k\geq 2$, we can write $m(s)=a(s)-b(s)$ for some nonnegative, nondecreasing functions $a(s)$, $b(s)$ on $[0,2\pi]$, satisfying
$$\int_0^{2\pi}da(s)\leq\frac{k}{2}+1,\;\;\;\int_0^{2\pi}db(s)\leq\frac{k}{2}-1.\eqno{(1.2)}$$
Furthermore let $P_k(\beta)$ be the class of analytic functions in $E$ which have
the representation
$$h(z)=\frac{1}{2}\int_0^{2\pi}\frac{1+(1-2\beta)ze^{-is}}{1-ze^{-is}}dm(s)\eqno{(1.3)}$$
where $m\in M_k$.\vskip 2mm

If we set $a(s)=\tfrac{k+2}{4}\mu_1(s)$ and $b(s)=\tfrac{k-2}{4}\mu_2(s)$ where $\int_0^{2\pi}d\mu_j(s)=1$, $j=1$, 2. Then from (1.2) and (1.3) we have
$$\aligned
h(z)
&=\frac{k+2}{4}\int_0^{2\pi}\frac{1+(1-2\beta)ze^{-is}}{1-ze^{-is}}d\mu_1(s)\\
&\qquad -\frac{k-2}{4}\int_0^{2\pi}\frac{1+(1-2\beta)ze^{-is}}{1-ze^{-is}}d\mu_2(s)\\
&\equiv\frac{k+2}{4}p(z)-\frac{k-2}{4}q(z)
\endaligned\eqno{(1.4)}$$
where $p$, $q\in P(\beta)$.\vskip 2mm

Pinchuk \cite{BP} defined $P_k$ and proved that\vskip 2mm

{\bf Lemma 1.1.} (\cite{BP}) {\em
All functions in $P_k$ have positive real parts for
$$|z|=r<\frac{k-\sqrt{k^2-4}}{2}.$$
Furthermore, there exist functions in $P_k$ which do not have positive real parts in any larger disk.}\vskip 2mm

If we  choose $h\in P_k(\beta)$ and set $h(z)=\beta+(1-\beta)h_1(z)$ where $h_1\in P_k$, we have the following\vskip 2mm

{\bf Lemma 1.2.} {\em
All functions in $P_k(\beta)$ have positive real parts for
$$|z|=r(k,\beta)<\left\{
\begin{array}{ll}\,
\frac{(1-\beta)k-\sqrt{(1-\beta)^2k^2-4(1-2\beta)}}{2(1-2\beta)}&\mbox{if $\beta\neq\frac{1}{2}$},\\
\frac{2}{k}&\mbox{if $\beta=\frac{1}{2}$.}
\end{array}
\right.\eqno{(1.5)}$$
The function
$$H(z)=\frac{k+2}{4}L_{0,\beta}(-z)-\frac{k-2}{4}L_{0,\beta}(z)\eqno{(1.6)}$$
shows that there exist functions in $P_k(\beta)$ which do not have positive real parts in any larger disk.}\vskip 2mm

In the next section we introduce certain iterations of $P_k(\beta)$ and with a brief discussion that will lead to new classes of analytic functions having boundary rotations bounded by $k\pi$.
\bigskip

{\bf 2. Functions of the classes $\Phi_{\sigma,n,k}^j(\beta)$, $j=1$, 2}
\medskip

In \cite{BO,KO}, the authors identified the following iterated integral transformation of functions in the class $P$.
\medskip

{\bf Definition 2.1.} Let $p\in P$. Let $\sigma>0$ be real number $n\geq 1$. Then, for $z\in E$ define
$$\phi_{\sigma,n}^j(p(z))=\int_0^z\lambda_{\sigma,n}^j(z,t)\phi_{\sigma,n-1}^j(p(t))dt,\;\;\;j=1,\;2\eqno{(2.1)}$$
where
$$\lambda_{\sigma,n}^1(z,t)=\frac{\sigma t^{\sigma-1}}{z^\sigma},$$
and
$$\lambda_{\sigma,n}^2(z,t)=\frac{(\sigma-(n-1))t^{\sigma-n}}{z^{\sigma-(n-1)}},\;\;\;\,\sigma-(n-1)>0$$
with $\phi_{\sigma,0}^j(p(z))=p(z)$. Changes in notations (cf. \cite{BO, KO}) became necessary only in order to unify our discussions of the two transformations. For any $p\in P$ defined by (1.1), it is known \cite{BO,KO} that the transformations $\phi_{\sigma,n}^j(p(z))$ have series representations
$$\phi_{\sigma,n}^j(p(z))=1+(1-\beta)\sum_{l=1}^\infty c_{l,n}^jz^l,\;\;\;j=1,\;2\eqno{(2.2)}$$
where
$$c_{l,n}^1=\left(\frac{\sigma}{\sigma+l}\right)^n$$
and
$$c_{l,n}^2=\frac{\sigma
(\sigma-1)...(\sigma-(n-1))}{(\sigma+l)(\sigma+l-1)...(\sigma+l-(n-1))}.$$
\vskip 2mm

Furthermore these transformations preserve many geometric structures of the family $P$; particularly the positivity of the real parts, compactness, convexity and subordination. They are iterative and are closely associated with certain families analytic and univalent functions involving the well known Salagean and Rucheweyh derivatives (see \cite{BO,KO}) and have been used effectively, and elegantly too, to characterize them. Also, in \cite{KB}, we have used the tranformation $\phi_{\sigma,n}^1(p(z))$ together with a method of Nehari and Netanyahu to determine the best possible coefficient bounds for some classes of functions. Further characterizations can be found in the articles.\vskip 2mm

Since for $0\leq\beta<1$, Re $p(z)>\beta$ implies Re $\phi_{\sigma,n}^j(p(z))>\beta$ (see \cite{BO,KO}), we would suppose $p\in P(\beta)$ and denote by $\Phi_{\sigma,n}^j(\beta)$ the classes of these transformations. If in (1.4) we replace $p(z)$ and $q(z)$ by their respective iterations (i.e. $\phi_{\sigma,n}^j(p(z))$ and $\phi_{\sigma,n}^j(q(z))$) we come to a more general class of functions $\Phi_{\sigma,n,k}^j(\beta)$ consisting of functions of the form:
$$\phi_{\sigma,n}^j(h(z))=\frac{k+2}{4}\phi_{\sigma,n}^j(p(z))-\frac{k-2}{4}\phi_{\sigma,n}^j(q(z)).\eqno{(2.3)}$$

Let us note the equivalent classes: 
$$\Phi_{\sigma,0,k}^j(\beta)\equiv P_k(\beta),\;\;\Phi_{\sigma,0,k}^j(0)\equiv P_k,$$
$$\Phi_{\sigma,0,2}^j(\beta)\equiv P(\beta),\;\text{and}\;\Phi_{\sigma,0,2}^j(0)\equiv P.$$

{\bf Theorem 2.2} {\em All functions in $\Phi_{\sigma,n,k}^j(\beta)$ have integral representation
$$\phi_{\sigma,n}^j(h(z))=\frac{1}{2}\int_0^{2\pi}\left[1+2(1-\beta)\sum_{l=1}^\infty c_{l,n}^jz^le^{-lis}\right]dm(s),\;\;\;j=1,2\eqno{(2.4)}$$
where $c_{l,n}^j$ are those defined for series $(2.2)$.}
\begin{proof}
From (1.3) and (2.1) we have
$$\aligned
\phi_{\sigma,n}^j(h(z))
&=\int_0^z\lambda_{\sigma,n}^j(z,t)\phi_{\sigma,n-1}^j(h(t))dt,\\
&=\int_0^z\lambda_{\sigma,n}^j(z,t)\phi_{\sigma,n-1}^j\left[\frac{1}{2}\int_0^{2\pi}\frac{1+(1-2\beta)te^{-is}}{1-te^{-is}}dm(s)\right]dt\\
&=\frac{1}{2}\int_0^{2\pi}\left[\int_0^z\lambda_{\sigma,n}^j(z,t)\phi_{\sigma,n-1}^j\left(\frac{1+(1-2\beta)te^{-is}}{1-te^{-is}}\right)dt\right]dm(s).
\endaligned\eqno{(2.5)}$$
If we write the kernel of the inner integral in series form and apply the transformations successively, we obtain the desired expression.
\end{proof}

{\bf Remark 2.3.} {\em The representations $(2.4)$ are the Herglotz's for functions in $\Phi_{\sigma,n,k}^j(\beta)$, $j=1$, 2.}\vskip 2mm

{\bf Theorem 2.4.} {\em For $h\in P(\beta)$, let $\phi_{\sigma,n}^j(h(z))\in\Phi_{\sigma,n,k}^j(\beta)$. Then we have the best possible lower bound
$$Re\;\phi_{\sigma,n}^j(h(z))\geq 1+(1-\beta)\sum_{l=1}^\infty\left(2c_{2l,n}^jr-kc_{2l-1,n}^j\right)r^{2l-1}.$$
where $c_{2l,n}^j$ and $c_{2l-1,n}^j$ are appropriately defined from those for series $(2.2)$.}
\begin{proof}
Let $z=re^{i\theta}$ and $t=\rho e^{i\theta}$ with $0<\rho<r<1$ in the second equation of (2.5), we have
$$\phi_{\sigma,n}^j(h(re^{i\theta}))
=\int_0^r\lambda_{\sigma,n}^j(r,\rho)\phi_{\sigma,n-1}^j(\Omega(\rho))d\rho$$
where
$$\aligned
\Omega(\rho)
&=\frac{1}{2}\int_0^{2\pi}\frac{1+(1-2\beta)\rho e^{i(\theta-s)}}{1-\rho e^{i(\theta-s)}}dm(s)\\
&=\frac{1}{2}\int_0^{2\pi}\left(\beta+(1-\beta)\frac{1+\rho e^{i(\theta-s)}}{1-\rho e^{i(\theta-s)}}\right)dm(s).
\endaligned$$
Thus we have
$$Re\;\phi_{\sigma,n}^j(h(re^{i\theta}))
=\int_0^r\lambda_{\sigma,n}^j(r,\rho)\phi_{\sigma,n-1}^j(Re\;\Omega(\rho))d\rho.\eqno{(2.6)}$$
Assume $\int_0^{2\pi}|dm(s)|=k$, then we have
$$\aligned
Re\;\Omega(\rho)
&=\frac{1}{2}\int_0^{2\pi}\left(\beta+(1-\beta)\frac{1-\rho^2}{1-2\rho\cos(\theta-s)+\rho^2}\right)[da(s)-db(s)]\\
&\geq\beta+(1-\beta)\left[\left(\frac{k}{2}+1\right)\frac{1-\rho^2}{1+2\rho+\rho^2}-\left(\frac{k}{2}-1\right)\frac{1-\rho^2}{1-2\rho+\rho^2}\right]\\
&=\beta+(1-\beta)\left(\frac{\rho^2-k\rho+1}{1-\rho^2}\right)\\
&=\frac{(1-2\beta)\rho^2-(1-\beta)k\rho+1}{1-\rho^2}.
\endaligned\eqno{(2.7)}$$
In series form we have
$$Re\;\Omega(\rho)\geq 1+(1-\beta)\sum_{l=1}^\infty\left(2\rho-k\right)\rho^{2l-1}.$$
Now applying (2.6) appropriately we have the results. And if $\int_0^{2\pi}|dm(s)|=k^\prime<k$, then we have Re $\Omega(\rho)\geq\tfrac{(1-2\beta)\rho^2-(1-\beta)k^\prime\rho+1}{1-\rho^2}>\tfrac{(1-2\beta)\rho^2-(1-\beta)k\rho+1}{1-\rho^2}$. The extremal functions are $\phi_{\sigma,n}^j(H(z))$, $H(z)$ given by (1.6).
\end{proof}

Two important corollaries ($j=1$, 2) which follow from the proof of the above theorem are results similar to Lemma 1.2, which are\vskip 2mm

{\bf Corollary 2.5.} {\em All functions in $\Phi_{\sigma,n,k}^j(\beta)$ have positive real parts for $|z|=r(k,\beta)$ given by $(1.5)$. This disk $|z|=r(k,\beta)$ is the largest possible.}
\begin{proof}
The proof easily follows from (2.7) and the functions $\phi_{\sigma,n,k}^j(H(z))$, $H(z)$ given by (1.6), show that there exist functions in $\Phi_{\sigma,n,k}^j(\beta)$ which do not have positive real parts in any larger disk.
\end{proof}

Following from (2.3) (and similar arguments as in \cite{BO,KO}) it is not difficult to see the inclusions\vskip 2mm

{\bf Theorem 2.6.} {\em 
$\Phi_{\sigma,n+1,k}^j(\beta)\subset\Phi_{\sigma,n,k}^j(\beta),\;\;\;\,n\in\mathbb{N}.$
}\vskip 2mm

{\bf Theorem 2.7.} {\em 
$\Phi_{\sigma,n,k}^j(\beta)\subset P_k(\beta),\;\;\;\,n\in\mathbb{N}.$
}\vskip 2mm

Finally, we make the following remarks.\vskip 2mm

{\bf Remark 2.8.} {\em 
For any real numbers $\sigma_1$, $\sigma_2>0$ and $n_1$, $n_2\in \mathbb{N}$ satisfying Definition $1.1$ and for $h\in P(\beta)$, we have
$$\phi_{\sigma_1,n_1}^j[\phi_{\sigma_2,n_2}^j(h(z)]=\phi_{\sigma_2,n_2}^j[\phi_{\sigma_1,n_1}^j(h(z)].$$
}\vskip 2mm

{\bf Remark 2.9.} {\em 
For any $h\in P(\beta)$, we have
$$\phi_{\sigma,n}^1(h(z)+\frac{z(\phi_{\sigma,n}^1(h(z))^\prime}{\sigma}=\phi_{\sigma,n-1}^1(h(z),$$
and
$$\phi_{\sigma,n}^2(h(z)+\frac{z(\phi_{\sigma,n}^2(h(z))^\prime}{\sigma-(n-1)}=\phi_{\sigma,n-1}^2(h(z).$$
}
\bigskip

{\bf 3. Functions of Bounded Boundary Rotation}
\medskip

Let $A$ be the class of functions
$$f(z)=z+a_2z^2+...$$
which are analytic in $E$. Throughout $\sigma$ and $n$ shall have their definitions as in Section 1. Then via the classes of analytic functions $P_k(\beta)$, we define the following classes of functions. 
\medskip

{\bf Definition 3.1.} {\em We say a function $f\in A$ is in the class $T_n^\sigma(k,\beta)$ if and only if
$$\frac{D^nf(z)^\sigma}{\sigma^nz^\sigma}\in P_k(\beta), \;\;\;z\in E.$$
}
\medskip

{\bf Definition 3.2.} {\em We say a function $f\in A$ is in the class $B_n^\sigma(k,\beta)$ if and only if
$$\frac{L_n^\sigma f(z)}{z}\in P_k(\beta), \;\;\;z\in E.$$}
\medskip

The operator $D^n:A\rightarrow A$ is the Salagean derivative operator defined as $D^nf(z)=D(D^{n-1}f(z))=z[D^{n-1}f(z)]^{\prime}$ with
$D^0f(z)=f(z)$. The index $\sigma>0$ means principal determinations only. Whereas, $L_n^\sigma:A\rightarrow A$ was defined in \cite{KO} (using the convolution $*$) as follows:
$$L_n^\sigma f(z)=(\tau_\sigma*\tau_{\sigma,n}^{(-1)}*f)(z).$$
where
$$\tau_{\sigma,n}(z)=\frac{z}{(1-z)^{\sigma-(n-1)}},\;\;\;\sigma-(n-1)>0,$$
$\tau_\sigma=\tau_{\sigma,0}$ and $\tau_{\sigma,n}^{(-1)}$ is such that
$$(\tau_{\sigma,n}*\tau_{\sigma,n}^{(-1)})(z)=\frac{z}{1-z}.$$

It is remarkable to mention that $L_n^n\equiv\mathcal{D}^n$ where $\mathcal{D}^n$ is the well known Rusheweyh derivative operator. $T_n^\sigma(2,\beta)$ and $B_n^\sigma(2,\beta)$ respectively coincide with the classes $T_n^\sigma(\beta)$ and $B_n^\sigma(\beta)$, which are generalizations of several other classes of functions \cite{SA,RM,TH,RS,KY}. Also $T_1^1(k,0)$ and $B_1^1(k,0)$ both coincide with $B_k$ studied in \cite{PB}. The classes $T_n^\sigma(k,\beta)$ and $B_n^\sigma(k,\beta)$ consist of analytic and univalent functions in $A$ having boundary rotations bounded by $k\pi$ and respectively generalize the classes $T_n^\sigma(\beta)$ and $B_n^\sigma(\beta)$ in the same manner $U_k$ and $V_k$ generalize the well known classes of starlike and convex functions in the open unit disk (see \cite{KI,BP,PB}).\vskip 2mm

In the sequel, we study these classes of functions. The conciseness of our proofs is due to the following lemmas, which relate functions of these classes with those of $\Phi_{\sigma,n,k}^j(\beta)$ of the last section. The lemmas are consequences of Lemmas 4.2 and 2 of \cite{BO, KO} respectively.\vskip 2mm

{\bf Lemma 3.3.} {\em Let $f\in A$. Then the following are equivalent:

{\rm(i)} $f\in T_n^\sigma(k,\beta)$,

{\rm(ii)} $\frac{D^nf(z)^\sigma}{\sigma^nz^\sigma}\in P_k(\beta)$,

{\rm(iii)} $f(z)^\sigma/z^\sigma\in \Phi_{\sigma,n,k}^1(\beta)$.}\vskip 2mm

{\bf Lemma 3.4.} {\em Let $f\in A$. Then the following are equivalent:

{\rm(i)} $f\in B_n^\sigma(k,\beta)$,

{\rm(ii)} $\frac{L_n^\sigma f(z)}{z}\in P_k(\beta)$,

{\rm(iii)} $f(z)/z\in\Phi_{\sigma,n,k}^2(\beta)$.}\vskip 2mm

Our first results are the inclusions.\vskip 2mm

{\bf Theorem 3.5.} {\em For $n\in\mathbb{N}$,
$$T_{n+1}^\sigma(k,\beta)\subset T_n^\sigma(k,\beta),\;\;B_{n+1}^\sigma(k,\beta)\subset B_n^\sigma(k,\beta).$$
}
\begin{proof}
Let $f\in T_{n+1}^\sigma(k,\beta)$. By Lemma 3, $f(z)^\sigma/z^\sigma\in\Phi_{\sigma,n+1,k}^1(\beta)$. Hence by Theorem 3, we have $f(z)^\sigma/z^\sigma\in\Phi_{\sigma,n,k}^1(\beta)$. This implies $f\in T_n^\sigma(k,\beta)$ as required. The proof of the second part is similar.
\end{proof}

It is known that for $n\geq 1$ the classes $T_n^\sigma(\beta)$ and $B_n^\sigma(\beta)$ both consist only of univalent functions in the open unit disk. Based on this, the next results show that $|z|=r(k,\beta)$ given by (1.5) is the radius of univalence for both $T_n^\sigma(k,\beta)$ and $B_n^\sigma(k,\beta)$ when $n\geq 1$.\vskip 2mm

{\bf Theorem 3.6.} {\em If $f\in T_n^\sigma(k,\beta)$, then
$$Re\;\frac{D^nf(z)^\sigma}{z^\sigma}>0\eqno{(3.1)}$$
for $|z|=r(k,\beta)$ given by $(1.5)$. For any larger disk, there exist functions in $T_n^\sigma(k,\beta)$ which do not satisfy $(3.1)$. Similarly, if $f\in B_n^\sigma(k,\beta)$, then
$$Re\;\frac{L_n^\sigma f(z)}{z}>0\eqno{(3.2)}$$
for $|z|=r(k,\beta)$ given by $(1.5)$. For any larger disk, there exist functions in $T_n^\sigma(k,\beta)$ which do not satisfy $(3.2)$.}
\begin{proof}
The proofs follow easily from Lemmas 2.2, 3.3, and 3.4. The functions $f(z)$ defined by $f(z)^\sigma/z^\sigma=\phi_{\sigma,n}^1(H(z))$ and $f(z)/z=\phi_{\sigma,n}^2(H(z))$, $H(z)$ given by (1.6), respectively show that the results are sharp.
\end{proof}

Define the integral
$$F(z)^\kappa=\frac{c+\kappa}{z^c}\int_0^zt^{c-1}f(t)^\kappa dt,\;\;\;c+\kappa>0.$$

For $\kappa=\sigma$ we consider the integral $F$ in $T_n^\sigma(k,\beta)$ while it is considered in $B_n^\sigma(k,\beta)$ for $\kappa=1$. Our results are the following\vskip 2mm

{\bf Theorem 3.7.} {\em Both classes $T_n^\sigma(k,\beta)$ and $B_n^\sigma(k,\beta)$ are preserved under $F$.}
\begin{proof}
Note that $\kappa=\sigma$ for $T_n^\sigma(k,\beta)$ and $\kappa=1$ for $B_n^\sigma(k,\beta)$. Suppose $f$ belong to either class. Then by Lemmas 3.3 and 3.4, $f(z)^\kappa/z^\kappa$ belongs to $\Phi_{\sigma,n,k}^j(\beta)$, $j=1$, 2; that is there exists $h\in P(\beta)$ such that $f(z)^\kappa/z^\kappa=\phi_{\sigma,n}^j(h(z))$. Hence
$$\aligned
F(z)^\kappa/z^\kappa
&=\frac{c+\kappa}{z^{c+\kappa}}\int_0^zt^{c+\kappa-1}\left(f(t)^\kappa/t^\kappa\right)dt,\;\;\;c+\kappa>0.\\
&=\frac{c+\kappa}{z^{c+\kappa}}\int_0^zt^{c+\kappa-1}\phi_{\sigma,n}^j(h(t))dt\\
&=\phi_{\gamma,1}^j[\phi_{\sigma,n}^j(h(t))],\;\;\;\gamma=c+\kappa\\
&=\phi_{\sigma,n}^j[\phi_{\gamma,1}^j(h(t))].
\endaligned$$
so that by Remark 2.8 and Lemmas 3.3 and 3.4 again we have the results.
\end{proof}

{\bf Theorem 3.8.} {\em Let $c_{l,n}^j$, $j=1$, $2$ have their definitions as in $(2.2)$. If $f\in T_n^\sigma(k,\beta)$, then
$$Re\;\frac{f(z)^\sigma}{z^\sigma}\geq 1+(1-\beta)\sum_{k=1}^\infty\left(2c_{2l,n}^1r-kc_{2l-1,n}^1\right)r^{2l-1}.$$
Similarly, if $f\in B_n^\sigma(k,\beta)$, we have
$$Re\;\frac{f(z)}{z}\geq 1+(1-\beta)\sum_{k=1}^\infty\left(2c_{2l,n}^2r-kc_{2l-1,n}^2\right)r^{2l-1}.$$
The bounds are the best possible.}
\begin{proof}
Take $f(z)^\sigma/z^\sigma=\phi_{\sigma,n}^1(h(z))$ and $f(z)/z=\phi_{\sigma,n}^2(h(z))$ in Theorem 2.6. The extremal functions are defined by $f(z)^\sigma/z^\sigma=\phi_{\sigma,n}^1(H(z))$ and $f(z)/z=\phi_{\sigma,n}^2(H(z))$, $H(z)$ given by (1.6).
\end{proof}

{\bf Theorem 3.9.} {\em Let $c_{l,n}^j$, $j=1$, $2$ have their definitions as in $(2.2)$. If $f\in T_n^\sigma(k,\beta)$, then
$$Re\;\frac{f(z)^{\sigma-1}f^\prime(z)}{z^{\sigma-1}}\geq 1+(1-\beta)\sum_{k=1}^\infty\left(2c_{2l,n-1}^1r-kc_{2l-1,n-1}^1\right)r^{2l-1}.$$
Similarly, if $f\in B_n^\sigma(k,\beta)$, we have
$$Re\;\frac{(\sigma-n)\frac{f(z)}{z}+f^\prime(z)}{\sigma-(n-1)}\geq 1+(1-\beta)\sum_{k=1}^\infty\left(2c_{2l,n-1}^2r-kc_{2l-1,n-1}^2\right)r^{2l-1}.$$
The bounds are the best possible.}
\begin{proof}
Take $f(z)^\sigma/z^\sigma=\phi_{\sigma,n}^1(h(z))$ and $f(z)/z=\phi_{\sigma,n}^2(h(z))$. On differentiation and using Remark 3, we have the results. The extremal functions are those defined by $f(z)^\sigma/z^\sigma=\phi_{\sigma,n}^1(H(z))$ and $f(z)/z=\phi_{\sigma,n}^2(H(z))$, $H(z)$ given by (1.6).
\end{proof}

{\bf Remark 3.10.} {\em Variants of our results for special cases of the parameters $\sigma$, $n$, $k$ and $\beta$ can be derived by specifying them.}

\medskip

{\em Acknowledgements.} This work was carried out at the Centre for Advanced Studies in Mathematics, CASM, Lahore University of Management Sciences, Lahore, Pakistan during the author's postdoctoral fellowship at the Centre. The author is indebted to all staff of CASM for their hospitality, most especially Prof. Ismat Beg. 
\bigskip

 \bibliographystyle{amsplain}

\end{document}